\documentclass[12pt]{article}
\usepackage{latexsym,amssymb}
\usepackage{amsmath}
\usepackage{color}
\usepackage{cite}

\newtheorem{Theorem}{\bf Theorem}[section]
\newtheorem{Lemma}{\bf Lemma}[section]
\newtheorem{Proposition}{\bf Proposition}[section]
\newtheorem{Corollary}{\bf Corollary}[section]
\newtheorem{Remark}{\bf Remark}[section]
\newtheorem{Example}{\bf Example}[section]
\newtheorem{Definition}{\bf Definition}[section]

\newenvironment{theorem}{\begin{Theorem}$\!\!\!$}{\end{Theorem}}
\newenvironment{lemma}{\begin{Lemma}$\!\!\!$}{\end{Lemma}}
\newenvironment{proposition}{\begin{Proposition}$\!\!\!$}{\end{Proposition}}

\newenvironment{remark}{\begin{Remark}$\!\!\!$}{\end{Remark}}

\newenvironment{definition}{\begin{Definition}$\!\!\!$}{\end{Definition}}

\numberwithin{equation}{section}

\def\XXint#1#2#3{{\setbox0=\hbox{$#1{#2#3}{\int}$}
\vcenter{\hbox{$#2#3$}}\kern-.5\wd0}}

\begin{document}

\title{Blow-up rate of sign-changing solutions \\to nonlinear parabolic systems}
\author{Erbol Zhanpeisov}
\date{}
\maketitle
\begin{abstract}
We present a  blow-up rate estimate for a solution to the parabolic Gross-Pitaevskii and related systems on entire space with Sobolev subcritical nonlinearity. We extend the results of [Y. Giga, S. Matsui and S. Sasayama, Indiana Univ. Math. J.  {53}  (2004), 483--514]  to the parabolic systems. 
\end{abstract}

\section{Introduction}
Let $U= (u_1,\dots,u_M)$ be a (classical) solution to the Cauchy problem for  a semilinear parabolic system
\begin{equation}
\tag{P}
\left\{
\begin{array}{ll}
\partial_t U - \Delta U = (\nabla G)(U) \quad & {\rm in} \quad {\bf R}^N \times(0,T) , \vspace{3pt}\\
U(x,0) = U_0 (x)& {\rm in}\quad {\bf R}^N,\vspace{3pt} 
\end{array}
\right.
\end{equation}
where $N\geq 1$, $M\geq 1$, $T>0$ and $U_0\in BC^2({\bf R}^N)$. Here  
\begin{equation}
\label{eq:1.1}
G(U)=\frac{1}{2(r+1)}\sum_{i,j=1}^M \beta_{ij}|u_i|^{r+1}|u_j|^{r+1} \quad \mbox{with} \quad r>0, 
\end{equation}
 where $\{\beta_{ij} \} \subset {\bf R}$ satisfies 
\begin{equation}
\label{eq:1.2}
\beta_{ij}\geq0,\quad \beta_{ii}>0, \quad \beta_{ij}=\beta_{ji}. 
\end{equation}
Then $\{u_i\}_{i=1} ^M$ satisfies 
\begin{align}
\label{eq:1.3}
& \partial_t u_i - \Delta u_i = \displaystyle\sum_{j=1}^{M}\beta _{ij}|u_i|^{r-1}|u_j|^{r+1}u_i  \quad {\rm in}\quad {\bf R}^N \times(0,T),\\ 
& u_i(x,0) = u_{i,0} (x) \quad {\rm in}\quad {\bf R}^N. \nonumber 
\end{align}
Parabolic system~\eqref{eq:1.3} can be regarded as a parabolic generalization of the Gross-Pitaevskii system
\begin{equation*}
\frac{1}{\sqrt{-1}} \partial_t u_i - \Delta u_i = \sum_{j=1}^M \beta_{ij} u_i |u_j|^2 \quad {\rm in} \quad {\bf R}^N \times(0,T), \quad {\rm for} \quad i=1, \dots, M,  
\end{equation*} 
and it is called parabolic Gross-Pitaevskii and related system. Parabolic system \eqref{eq:1.3} has nice mathematical structure such as scaling invariance and energy structure. In this paper, thanks to these nice properties of parabolic system \eqref{eq:1.3}, we study a blow-up rate estimate for a solution to problem~(P). 

Blow-up rate of nonlinear parabolic problems  has been studied in many papers, see e.g., \cite{C, CMR1, CMR2, CMR3, CRS, DMW, FS, FHV, FM, FI1, FI2, FI3, FIM, GK1, GK2, GK3, GMS1, GMS2, HV4, IM, M, MM1, MM2, MM3, MS, MRS, MZ2, MZ3, M1, M2, M3, M4, M5, PQS2, Q1, Q2, QS2, Sc, Se}. 
(See also the monograph \cite{QS}.)
Let us consider  problem~(P) with $M=1$, that is  \begin{equation}
\label{eq:1.4}
\left\{
\begin{array}{ll}
u_t-\Delta u =|u|^{p-1} u &  {\rm in}\quad {\bf R}^N\times(0,T) ,\vspace{3pt}\\
u(x,0) = u_{0} (x)& {\rm in}\quad {\bf R}^N, 
\end{array}
\right.
\end{equation}
where $p>1$. 
Assume that
\begin{equation*}
\limsup_{t \to T}\,  \|u(t)\|_{L^\infty ({\bf R}^N)} = \infty.
\end{equation*}
We say that the blow-up is of type I if 
\begin{equation*}
\limsup_{t \to T}\, (T-t)^{\frac{1}{p-1}}\|u(t)\|_{L^\infty ({\bf R}^N)}<\infty, 
\end{equation*} 
and of type II otherwise. 
It is known that if $1<p<p_S$, then the solution $u$ to problem~\eqref{eq:1.4} exhibits type I blow-up, where 
$$
p_S:=\frac{N+2}{N-2}\quad\mbox{if}\quad N\ge 3,
\qquad
p_S:=\infty\quad\mbox{if}\quad N=1,2.
$$
See \cite{GK2, GMS1, GMS2}. 
On the other hand, if $p\geq p_S$, the solution $u$ to problem \eqref{eq:1.4} does not necessarily exhibit type I blow-up. See e.g., \cite{C, CMR3, DMW, FHV, HV4, MM3, M4, Sc, Se}.

For the parabolic system, blow-up estimates are much less understood. See e.g., \cite{MZ3, P, Q3, PS, Q4}. Among others, Phan and Souplet \cite{PS} established a parabolic Liouville-type theorem 
for problem~(P) and proved that  
 a solution $U=(u_1,\dots,u_M)$ of problem~(P) satisfies the type~I blow-up estimate 
\begin{equation}
\label{eq:1.5}
\|U\|_{L^\infty({\bf R}^N)}\leq 
C(T-t)^{-\frac{1}{p-1}},\quad
0<t<T,
\end{equation}
provided that $1<p:=2r+1<p_B$ and 
$U\ge 0$ in ${\bf R}^N\times(0,T)$, that is, 
$u_i\ge 0$ in ${\bf R}^N\times(0,T)$, where $i=1,\dots,M$.
Here
$$
p_B:=\frac{N(N+2)}{(N-1)^2}\quad\mbox{if}\quad N\ge 3,
\qquad
p_B:=\infty\quad\mbox{if}\quad N=1,2.
$$
On the other hand, 
it seems difficult to apply their arguments to problem~(P) 
without the restriction of the sign of $U$, 
since the parabolic Liouville-type theorem acts only on nonnegative solutions. 
In this paper we develop the arguments in \cite{GMS1, GMS2} and  obtain the type I blow-up estimate~\eqref{eq:1.5} without the restriction of the sign of $U$ for $1<p<p_S$. 
Note that $p_S>p_B$ if $N\ge 3$. 

We formulate a solution to problem~(P).
\begin{definition}
\label{Definition:1.1}
Let $M=1,2,\dots$, $T>0$ and $U_0\in BC^2({\bf R}^N)$. 
A function~$U=(u_1,\dots,u_M)$ in ${\bf R}^N\times[0,T)$ 
is said a solution to problem~{\rm (P)} if 
\begin{equation*}
\text{$u_i$, $\partial_t u_i$, $\nabla u_i$, $\nabla ^2 u_i$ are bounded
and continuous on ${\bf R}^N \times [0,\tau]$}
\end{equation*}
for $0<\tau<T$ and $i=1,\dots,M$ and 
$U=(u_1,\dots,u_M)$ solves {\rm (P)}. 
\end{definition}

We are ready to state the main result of this paper. 
\begin{theorem}
\label{Theorem:1.1}
Let $M=1,2,\dots$, $T>0$ and $U_0\in BC^2({\bf R}^N)$. 
Assume \eqref{eq:1.1} and \eqref{eq:1.2}. Let $U$ be a solution to problem~{\rm (P)}.
Then there exists  $C>0$ such that 
\begin{equation*}
\| U(t)\|_{L^\infty ({\bf R}^N)} \leq C(T-t)^{-\frac {1}{p-1}},
\quad 0<t<T, 
\end{equation*}
provided that $1<p:=2r+1<p_S$. 
\end{theorem}

We explain the idea of the proof of Theorem~\ref{Theorem:1.1}. Let $U$ be a solution in ${\bf R}^N \times [0,T)$. We may assume $T=1$ since problem~(P) possesses a similarity transformation such as 
\begin{equation*}
Z(x,t) = T^{\frac{1}{p-1}} U(\sqrt T x, Tt).
\end{equation*} 
As in \cite{GK2}, for any $a\in {\bf R}^N $, we introduce the following rescaled function of $U$ 
\begin{equation*}
W^a (y,s) := (T-t)^\beta U(x,t)
\end{equation*}
with
\begin{equation*}
y= \frac{x-a}{\sqrt{T-t}},\quad s=-\log (T-t),\quad \beta = \frac{1}{p-1}.
\end{equation*}
Then  $W^a(y,s)=(w_1,\dots,w_M)$ satisfies  
\begin{equation}
\label{eq:1.6}
\partial _s w_i - \Delta w_i+ \frac{1}{2}y\cdot \nabla w_i +\beta w_i - \sum_{j=1} ^M \beta_{ij} |w_i|^{r-1} |w_j|^{r+1}w_i =0, 
\end{equation} 
that is
\begin{equation}
\label{eq:1.7}
\rho \partial_s w_i - \nabla \cdot (\rho \nabla w_i) + \beta \rho w_i - \rho (\nabla G)(W^a) = 0  \quad {\rm in}\,\,\, {\bf R}^N\times(0,\infty)
\end{equation}
for $i=1,\dots, M$, where $\rho(y)=\exp(-|y|^2/4)$.
Furthermore, Definition~\ref{Definition:1.1} implies that 
\begin{equation}
\label{eq:1.8}
\begin{split}
&\text{$w_i$, $(1+ |y|)^{-1} \partial_s w_i$, $\nabla w_i$, $\nabla ^2 w_i$ are bounded
and continuous }\\
&\text{on ${\bf R}^N \times [0,S]$ for $i=1,2,\dots,M$ and  any $S<\infty$.}
\end{split}
\end{equation}
We introduce the global energy of $W$
\begin{equation}
\label{eq:1.9}
E[W](s) := \frac{1}{2} \int _{{\bf R}^N} \left(|\nabla W|^2 + \beta  |W|^2\right)\rho \, dy -  \int _{{\bf R}^N} G(W)\rho \,dy
\end{equation}
and the local energy of $W$
\begin{equation}
\label{eq:1.10}
E_{\varphi}[W] := \frac{1}{2}\int_{{\bf R}^N } \varphi^2 \left(|\nabla W|^2 + \beta |W|^2\right)\rho \,dy -\int_{{\bf R}^N } \varphi^2 G(W)\rho \,dy, 
\end{equation}
where 
$$
|\nabla W|^2 =\displaystyle{\sum _{i=1}^{M}} |\nabla w_i|^2
\quad\mbox{and}\quad
|W|^2 = \displaystyle{\sum _{i=1}^{M}} | w_i|^2.
$$
Here $\varphi(y)$ is a cut-off function.
Developing the argument in \cite{GMS1}, we prove  the boundedness of the global energy and the local energy of $W$, and obtain
\begin{equation}
\label{eq:1.11}
\sup _{s\geq0}\, \int _s^{s+1} \left(\int_{B(R)}|W|^{p+1}\rho \,dy \right)^q d\tau \leq C_{q,R}
\end{equation}
for any $R>0$ and $q\geq 2$. Here  $C_{q,R}$ is a constant  depending on $q$ and $R$. In \cite{GMS1}, for the case of $M=1$, the authors combined \eqref{eq:1.11} and the interior regularity theorem for parabolic equations to obtain the $L^\infty_{loc}$ estimate of $W$.  
However, this argument fails for  parabolic system \eqref{eq:1.3}  in the case where $M\geq2$ and $0<r<1$ because of the nonlinearity of $|w_i|^{r-1}$. In this paper, in order to  overcome the difficulty, we introduce a function $w$  defined by
\begin{equation*}
w:= \sum_{i=1}^M |w_i|, 
\end{equation*}
which is a weak subsolution to 
\begin{equation*}
\partial _s w \leq \Delta w- \frac{1}{2}y\cdot \nabla w -\beta w + \sum_{i,j=1} ^M \beta_{ij} w^{2r+1}\quad \mbox{in} \quad {\bf R}^N\times (0,\infty).
\end{equation*}
 Then we apply the parabolic regularity theorems to the weak subsolution $w$  and  obtain the estimate 
\begin{equation*}
\|w\|_{L^\infty (B_R)}\leq C \quad {\rm for } \quad s\geq 0.
\end{equation*}
This enables us to obtain the type I estimate, and  Theorem \ref{Theorem:1.1} follows.  

The rest of this paper is organized as follows. In Sections 2 and 3 we obtain the estimate for the global energy and the local energy respectively. In Section 4 we prove \eqref{eq:1.11} and Theorem~\ref{Theorem:1.1}.


\section{ Estimate of the global energy}
In this section we obtain the monotonicity of the energy $E[W](s)$ and its related inequalities.  
In what follows we write $F(U):=(\nabla G)(U)$ and often use 
the following structure conditions on $G$: 
\begin{equation}
\label{eq:2.1}
\left\{
\begin{array}{ll}
G(\lambda U) = \lambda^{p+1}G(U)\quad & \mbox{for $\lambda \ge 0$ and $U \in {\bf R}^M$};\vspace{3pt}\\
G(U)>0\quad & \mbox{for $U\in {\bf R}^M \setminus \{0\}$}. 
\end{array}
\right. 
\end{equation}
It follows from \eqref{eq:2.1} that 
\begin{equation}
\label{eq:2.2}
\left\{
\begin{array}{ll}
F(\lambda U) = \lambda^{p} F(U)\quad & \mbox{for $\lambda \ge 0$ and $U \in {\bf R}^M$};\vspace{3pt}\\
|F(U)| \leq C_F |U|^p \quad & \mbox{for $U\in {\bf R}^M$};\vspace{3pt}\\
 c_G |U|^{p+1}\leq G(U) \leq C_F |U|^{p+1} \quad & 
 \mbox{for $U\in {\bf R}^M$};\vspace{3pt}\\
U \cdot F(U) = (p+1)G(U)\quad & \mbox{for $U\in {\bf R}^M$}. 
\end{array}
\right.
\end{equation}
Here  $c_G$ and  $C_F $ are positive constants such that  
$$
c_G=\min _{|U|=1} G(U)\quad\mbox{and}\quad C_F=\max _{|U|=1} |F(U)|.
$$

In Propositions~\ref{Proposition:2.1} and \ref{Proposition:2.2}
we show the monotonicity and the nonnegativity of the global energy $E[W](s)$, respectively. 
\begin{proposition}
\label{Proposition:2.1}
Let $W$ satisfy \eqref{eq:1.6} with \eqref{eq:1.8}. 
Then 
\begin{align}
\label{eq:2.3}
& \frac{1}{2} \frac{d}{ds} \int _{{\bf R}^N} |W|^2 \rho \,dy = -2 E[W](s) + (p-1)\int _{{\bf R}^N} G(W)\rho \,dy,\\ 
\label{eq:2.4}
& \frac{d}{ds} E[W](s) = - \int_{{\bf R}^N} |W_s|^2 \rho \,dy, 
\end{align}
for $s\geq0$. 
\end{proposition}
{\bf Proof.}
It follows from \eqref{eq:1.7}, \eqref{eq:1.9} and \eqref{eq:2.2} that
\begin{equation}
\label{eq:2.5}
\begin{split}
\frac{1}{2} \frac{d}{ds} \int _{{\bf R}^N} |W|^2 \rho \,dy 
& = \int _{{\bf R}^N} W\cdot W_s \rho \,dy \\
& = -2 E[W](s) + (p-1)\int _{{\bf R}^N} G(W)\rho \,dy.
\end{split}
\end{equation}
This implies \eqref{eq:2.3}. 
Furthermore, since  
\begin{align*}
&  \frac{d}{ds} \int _{{\bf R}^N} G(W) \rho \,dy = \int _{{\bf R}^N} F(W)\cdot W_s \rho \,dy ,\\
& \frac{1}{2}\frac{d}{ds} \int _{{\bf R}^N} |\nabla W|^2 \rho \,dy
 = -\sum _{i=1}^M\int _{{\bf R}^N} w_{is} \nabla \cdot (\rho \nabla w_i)  \,dy,
\end{align*}
by \eqref{eq:2.5} we have \eqref{eq:2.4}. 
The proof is complete. 
$\Box$
\begin{proposition}
\label{Proposition:2.2}
Let $W$ satisfy \eqref{eq:1.6} with \eqref{eq:1.8}. 
Then 
\begin{equation}
\label{eq:2.6}
0\leq E[W](s) \leq E[W](0)\quad\mbox{for}\quad s\ge 0.
\end{equation} 
\end{proposition}
{\bf Proof.}
Due to \eqref{eq:2.4} it suffices to prove the nonnegativity  of $E[W](s)$. 
By \eqref{eq:2.2}, 
applying Jensen's inequality to \eqref{eq:2.5}, 
we have
\begin{equation}
\label{eq:2.7}
\begin{split}
\frac{1}{2} \frac{d}{ds} \int _{{\bf R}^N} 
|W(s)|^2 \rho\,dy \geq -2E[W](s) 
+ c_1 \left( \int _{{\bf R}^N} |W(s)|^2 \rho\,dy\right)^{\frac{p+1}{2}},
\end{split}
\end{equation}
where $c_1$ is a positive constant depending on $c_G$. 
If $E[W](s_*)<0$ for some $s_*\ge 0$, 
then by \eqref{eq:2.7} and the monotonicity of $E[W](s)$ 
we find $s_*'\in(s_*,\infty)$ such that  
$$
\int _{{\bf R}^N}|W(s)|^2 \rho\,dy\to\infty
\quad\mbox{as}\quad s\to s_*'.
$$
This contradicts the global existence of $W$. 
Thus Proposition~\ref{Proposition:2.2} follows.
$\Box$\vspace{5pt}

\begin{remark}
\label{Remark:2.1}
$E[W](0)$ is uniformly bounded with respect to $a$ since $U_0\in BC^2({\bf R}^N)$. Indeed, we have 
\begin{equation*}
M_0 := \sup_{a\in {\bf R}^N} ~E[W](0)\leq C \sup_{x\in{\bf R}^N}~ (|U_0(x)|^2  + |\nabla U_0(x)|^2 )< \infty.
\end{equation*}
\end{remark}

For any Banach space $X$ and $f\in X$, 
we denote by $\|f;X\|$ 
the norm of $f$ in $X$.
For any domain $\Omega\subset{\bf R}^N$ and $1\le q\le\infty$, 
we write $L^q_\rho(\Omega)=L^q(\Omega,\rho\,dy)$.
In particular, 
\begin{equation*}
 \|f;L_\rho ^2 (\Omega)\|:=  
 \left(\int_\Omega |f|^2\rho\,dy\right)^{\frac{1}{2}},
 \quad f\in L^2_\rho(\Omega).
\end{equation*}
Furthermore, we write 
$W^{1,2}_\rho(\Omega)=W^{1,2}(\Omega,\rho\,dy)$ 
and set 
\begin{equation*}
 \|f;W^{1,2}_\rho (\Omega)\|:=  
 \left(\int_\Omega (|\nabla f|^2+ \beta |f|^2)\rho \, dy\right)^{\frac{1}{2}},\quad f\in W^{1,2}_\rho(\Omega).
\end{equation*}
\begin{proposition}
\label{Proposition:2.3}
Let $W$ satisfy \eqref{eq:1.6} with \eqref{eq:1.8}. Then
\begin{equation}
\label{eq:2.8}
\int _0^\infty \|W_s;L_\rho ^2 ({\bf R}^N)\|^2 \,d\tau \leq E[W](0). 
\end{equation}
Moreover, there exist positive constants $K_1$, $K_2$ and $K_3$ depending only on $N$, $p$, $c_G$ and $M_0$ such that
\begin{align}
\label{eq:2.9}
& \sup_{s\geq 0}\, \|W;L_\rho ^2 ({\bf R}^N)\| \leq K_1,\\ 
\label{eq:2.10}
& \sup_{s\geq0}\, \int _s^{s+1} \|W;L_\rho ^{p+1} ({\bf R}^N)\|^{2(p+1)}\,d\tau \leq K_2,\\ 
\label{eq:2.11}
& \|W(s);W_\rho^{1, 2}({\bf R}^N)\|^2 \leq K_3 (1+\|W_s;L_\rho ^2 ({\bf R}^N)\|)\quad \mbox{for} \quad s\geq 0.
\end{align}
\end{proposition}
{\bf Proof.}
Inequality \eqref{eq:2.8} follows from  \eqref{eq:2.4} and \eqref{eq:2.6}.
We prove inequality \eqref{eq:2.9}. Let $f(s):=\|W;L^2_\rho ({\bf R}^N)\|$. By the Sobolev inequality we find $C> 0$ such that
\begin{equation}
\label{eq:2.12}
\sup_{[s,s+1]}\, f\leq C\left(\| f'; {L^2 (s,s+1)}\| + \| f;{L^2(s,s+1)}\| \right)^\theta \| f;{L^{2p}(s,s+1)}\| ^{1-\theta} 
\end{equation}
for $s\ge 0$, where $\theta=1/(p+1)$. By the definition of $f$ and Schwartz inequality  we have 
\begin{equation*}
|f'(s)|\leq \|W_s; L_\rho ^2 ({\bf R}^N)\|. 
\end{equation*}
This together with  \eqref{eq:2.8} implies 
\begin{equation}
\label{eq:2.13}
\| f';{L^2 (s,s+1)}\| \leq E[W](0)^{\frac{1}{2}}.
\end{equation}
We see  by \eqref{eq:2.6} and \eqref{eq:2.7} that
\begin{equation*}
2E[W](0) + f(s)f'(s) \geq c_1 f(s)^{p+1},  
\end{equation*}
where $c_1$ is as in \eqref{eq:2.7}.
If $f(s)>1$, then by Scwartz inequality we have 
\begin{equation*}
f(s)^{2p} \leq \frac{1}{c_1 ^2}\left(2\left|f'(s)\right|^2 + 8E[W](0)^2\right). 
\end{equation*}
Together with the case $f(s)\leq 1$ we get 
\begin{equation*}
\int _s^{s+1} f(\tau)^{2p} \,d\tau \leq 1+ 8c_1^{-2}E[W](0)^2 + 2c_1^{-2} \int_s^{s+1} \left|f'(\tau)\right|^2 \,d\tau. 
\end{equation*}
Thus by \eqref{eq:2.13} we have
\begin{equation}
\label{eq:2.14}
\|f ;{L^{2p}(s,s+1)}\| \leq \left(1+ 8c_1^{-2}E[W](0)^2 + 2c_1^{-2}E[W](0)\right)^{\frac{1}{2p}}. 
\end{equation}
Since $\| f;{L^2(s,s+1)}\| \leq \| f;{L^{2p} (s,s+1)}\|$, inequality \eqref{eq:2.9} follows from  \eqref{eq:2.12}, \eqref{eq:2.13} and \eqref{eq:2.14} with $K_1$ depending on $N$, $p$, $c_G$ and $M_0$.

\noindent By \eqref{eq:2.2} and \eqref{eq:2.5} we obtain
\begin{equation*}
\begin{split}
 c_G ^2 (p-1)^2 \left(\int _{{\bf R}^N} |W|^{p+1} \rho \, dy\right)^2 
& \leq \left(f(s)\|W_s; L_\rho ^2 ({\bf R}^N)\|+ 2 E[W](0)\right)^2\\
& \leq 2f(s)^2 \|W_s; L_\rho ^2 ({\bf R}^N)\|^2 + 8 E[W](0)^2.
\end{split}
\end{equation*}
Now \eqref{eq:2.10} follows from \eqref{eq:2.9} and \eqref{eq:2.10}.  

\noindent Finally we show \eqref{eq:2.11}. By the definition of $E[W]$,  \eqref{eq:2.3}, \eqref{eq:2.6} and  \eqref{eq:2.9} we get 
\begin{equation*}
\begin{split}
\|W(s);W_\rho^{1,2}({\bf R}^N)\|^2
& = 2E[W](s)+2\int_{{\bf R}^N}G(W)\rho \,dy \\
& =\frac{2}{p-1}\left( (p+1)E[W](s) + \int_{{\bf R}^N}W\cdot W_s\rho\, dy  \right)\\
& \leq \frac{2}{p-1}\left( (p+1)E[W](0) + K_1 \|W_s;L_\rho^2({\bf R}^N)\| \right).
\end{split}
\end{equation*}
Therefore inequality \eqref{eq:2.11} follows. 
$\Box$

\section{Estimates of local energy}
The aim of this section is  to prove the following two propositions. 
\begin{proposition}
\label{Proposition:3.1}
Let $W$ satisfy \eqref{eq:1.6} with \eqref{eq:1.8} and  $\psi \in BC^1({\bf R}^N)$. Then there exists  $L_1 >0$ such that
\begin{equation}
\label{eq:3.1}
E_\psi [W](s)\leq L_1, \quad \mbox{for}\quad s\ge 0.
\end{equation}
The constant $L_1$ depends only on $N$, $p$, $c_G$, $\| \psi \|_{BC^1({\bf R}^N)}$ and $M_0$. 
\end{proposition}

\begin{proposition}
\label{Proposition:3.2}
Let $W$ satisfy \eqref{eq:1.6} with \eqref{eq:1.8} and  $\psi \in BC^2({\bf R}^N)$ with  ${\rm supp}\  \psi \subset {B_R}$. Then there exists  $L_2>0$ such that
\begin{equation}
\label{eq:3.2}
E_\psi[W](s)\geq -L_2, \quad \mbox{for}\quad s\ge 0.
\end{equation}
The constant $L_2$ depends only on $N$, $p$, $c_G$, $R$, $\| \psi \|_{BC^2({\bf R}^N)}$ and $M_0$.
\end{proposition}

For this aim  we prepare the following three lemmas.

\begin{lemma}
\label{Lemma:3.1}
Let $W$ satisfy \eqref{eq:1.6} with \eqref{eq:1.8} and $\psi \in BC^1({\bf R}^N)$. Then 
\begin{align}
\label{eq:3.3}
& \frac{1}{2}\frac{d}{ds}  \int_{{\bf R}^N}  \psi^2 |W|^2 \rho\, dy \\
& = -2E_\psi [W](s) + (p-1)\int_{{\bf R}^N} \psi^2G(W)\rho \, dy -2\sum_{i=1}^M \int_{{\bf R}^N} \psi w_i \nabla \psi \cdot \nabla w_i\rho\, dy, \nonumber \\ 
\label{eq:3.4}
& \frac{d}{ds} E_\psi[W](s) = - \int_{{\bf R}^N} \psi^2|W_s|^2 \rho \,dy - 2\sum_{i=1}^M\int_{{\bf R}^N} \psi w_{is}\rho (\nabla w_i \cdot \nabla \psi)\,dy,
\end{align}
for $s\ge 0$. 
\end{lemma}
{\bf Proof.}
As in Proposition 2.1, we have  
\begin{equation}
\label{eq:3.5}
\begin{split}
\frac{1}{2} \frac{d}{ds} \int _{{\bf R}^N} \psi^2 |W|^2 \rho \,dy 
& = \int _{{\bf R}^N}\psi^2 W\cdot W_s \rho \,dy \\
& =\sum_{i=1}^M \int_{{\bf R}^N} \psi^2 w_i(\nabla\cdot (\rho \nabla w_i) -\beta \rho w_i + F_i (W)\rho) \,dy,
\end{split}
\end{equation} 
\begin{equation}
\label{eq:3.6}
\begin{split}
\frac{d}{ds} \int_{{\bf R}^N} & \left[\frac{1}{2} \beta \psi^2|W|^2 \rho  
 -\psi^2G(W)\rho \right]\, dy \\
& = \int_{{\bf R}^N} [\beta \psi^2 W\cdot W_s \rho -\psi^2 F(W)\cdot W_s \rho ]\,dy, 
\end{split}
\end{equation}
\begin{equation}
\label{eq:3.7}
\begin{split}
\frac{1}{2}\frac{d}{ds} \int_{{\bf R}^N}& \psi^2 |\nabla W|^2 \rho \,dy 
 = \sum_{i=1}^{M} \int_{{\bf R}^N} \rho \psi^2 \nabla w_i\cdot \nabla w_{is}\,dy \\
& = -\sum_{i=1}^M\int_{{\bf R}^N} \psi^2 w_{is}\nabla\cdot (\rho  \nabla w_i) \,dy - 2\sum_{i=1}^M\int_{{\bf R}^N} \psi w_{is} (\nabla w_i \cdot \nabla \psi)\rho\,dy. 
\end{split}
\end{equation}
Identity \eqref{eq:3.3} follows from \eqref{eq:3.5} and integration by parts. Furthermore, by \eqref{eq:3.6}, \eqref{eq:3.7} and \eqref{eq:1.7} we obtain \eqref{eq:3.4}. The proof is complete. 
$\Box$

\begin{lemma}
\label{Lemma:3.2}
Let $W$ satisfy \eqref{eq:1.6} with \eqref{eq:1.8} and  $\psi \in BC^1({\bf R}^N)$. Then there exist $L_3>0$ and $L_4>0$ such that 
\begin{align}
\label{eq:3.8}
& \frac{d}{ds} E_\psi [W](s) \leq L_3(1+ \| W_s(s); L_\rho ^2 ({\bf R}^N)\|),\\ 
\label{eq:3.9}
& \int_s^{s+1} E_\psi[W](\tau) \,d\tau \leq L_4,
\end{align} 
for $s\geq 0$. The constants $L_3$ and $L_4$ depend only on $N$, $p$, $c_G$, $\| \psi \|_{BC^1({\bf R}^N)}$ and $M_0$.
\end{lemma}
{\bf Proof.}
By \eqref{eq:3.4}, Cauchy's inequality and \eqref{eq:2.11}  we see that  
\begin{equation*}
\begin{split}
\frac{d}{ds} E_{\psi}[W](s)
& \leq -\int_{{\bf R}^N} \psi^2 |W_s|^2 \rho\, dy 
 +2\sum_{i=1}^M \int_{{\bf R}^N}|\psi||\nabla \psi|  |w_{is}| |\nabla w_i|\rho \,dy \\
& \leq -\frac{1}{2}\int_{{\bf R}^N} \psi^2 |W_s|^2 \rho \,dy 
+2\int_{{\bf R}^N}|\nabla \psi|^2 |\nabla W|^2\rho \,dy\\
& \leq 2K_3 \|\nabla \psi\|_\infty ^2 (1+\| W_s(s); L_\rho ^2 ({\bf R}^N)\|).
\end{split}
\end{equation*} 
This implies inequality \eqref{eq:3.8}.
By \eqref{eq:2.8},  \eqref{eq:2.11} and  Jensen's inequality we have 
\begin{equation*}
\begin{split}
\int_s^{s+1} E_\psi [W](\tau) \,d\tau 
& \leq \frac{\|\psi \|_{\infty} ^2}{2} \int_s^{s+1} \int_{{\bf R}^N} (|\nabla W|^2 + \beta |W|^2)\rho \,dy d\tau \\
& \leq K_3\|\psi \|_{\infty} ^2\left(1+ \int_s^{s+1} \|W_s(\tau);L_\rho^2 ({\bf R}^N)\| \,d\tau\right)\\
& \leq K_3\|\psi \|_{\infty} ^2 (1+ E[W](0)^{1/2}).
\end{split}
\end{equation*}
Therefore \eqref{eq:3.9} follows. The proof is complete.
$\Box\ $\vspace{5pt}

\begin{lemma}
\label{Lemma:3.3}
Let $W$ satisfy \eqref{eq:1.6} with \eqref{eq:1.8} and $\psi \in BC^2({\bf R}^N)$ with  ${\rm supp}\  \psi \subset {B_R}$. Then there exists   $L_5>0$ such that 
\begin{equation}
\label{eq:3.10}
\frac{1}{2}\frac{d}{ds}\int_{{\bf R}^N}  \psi^2 |W|^2 \rho\, dy  \geq -2E_\psi [W](s) + (p-1)\int_{{\bf R}^N} \psi^2G(W)\rho \,dy -L_5 .
\end{equation}
The constant $L_5$ depends only on $N$, $p$, $c_G$, $R$, $M_0$, $\| \psi \|_{BC^2({\bf R}^N)}$.
\end{lemma}
{\bf Proof.}
By \eqref{eq:3.3} it suffices to  prove 
\begin{equation*}
\sum_{i=1}^M\int_{{\bf R}^N} \psi w_{i}\rho (\nabla w_i \cdot \nabla \psi) \,dy\leq L_5.
\end{equation*} 
We get from integration by parts
\begin{equation*}
\begin{split}
\sum_{i=1}^M\int_{{\bf R}^N} \psi w_{i}\rho \nabla w_i \cdot \nabla \psi\,dy
& = - \sum_{i=1} ^M\int_{{\bf R}^N}  w_i \nabla \cdot (\psi w_i \rho \nabla \psi )\,dy \\
& = -\int_{{\bf R}^N} |\nabla \psi|^2 |W|^2 \rho \,dy -\sum_{i=1}^M\int_{{\bf R}^N} \psi w_i\rho \nabla w_i \cdot \nabla \psi \,dy \\
& - \int_{{\bf R}^N} \psi |W|^2 \Delta \psi \rho \,dy + \frac{1}{2}\int_{{\bf R}^N} \psi |W|^2 \rho (y \cdot \nabla \psi ) \,dy  . 
\end{split}
\end{equation*}
Therefore we obtain by \eqref{eq:2.9}
\begin{equation*}
\begin{split}
\sum_{i=1}^M\int_{{\bf R}^N} \psi w_{i}\rho \nabla w_i \cdot \nabla \psi\,dy
& = -\frac{1}{2} \int_{{\bf R}^N}  |\nabla \psi|^2 |W|^2 \rho \,dy -\frac{1}{2}\int_{{\bf R}^N} \psi \Delta \psi |W|^2 \rho \,dy \\
& + \frac{1}{4} \int_{{\bf R}^N} \psi |W|^2 \rho (y\cdot \nabla \psi ) \,dy\\
& \leq  \frac{1}{2} \|\psi\|_\infty \|\Delta \psi\|_\infty K_1^2 + \frac{1}{4} \|\psi\|_\infty \|\nabla \psi\|_\infty R K_1^2=:L_5.
\end{split}
\end{equation*} 
This implies \eqref{eq:3.10}.
$\Box$

We are ready to prove Propositions~\ref{Proposition:3.1} and \ref{Proposition:3.2}. \vspace{5pt}

\noindent{\bf Proof of Proposition \ref{Proposition:3.1}.}
Let $s\geq 1$. We set $s_1\in (s-1,s)$ such that 
\begin{equation}
\label{eq:3.11}
E_\psi[W](s_1) = \int_{s-1}^s  E_\psi[W](\tau)\,d\tau.
\end{equation}
Then we use \eqref{eq:3.8}, Jensen's inequality and \eqref{eq:2.8} to get
\begin{equation*}
\begin{split}
E_\psi[W](s)-E_\psi[W](s_1)  
& \leq L_3\left(1+ \int_{s-1}^s\| W_s(s); L_\rho ^2 ({\bf R}^N)\|\,d\tau \right)\\
& \leq L_3\left(1+ \left(E[W](0)\right)^{1/2}\right).
\end{split}
\end{equation*}
Combining this inequality with \eqref{eq:3.9} and \eqref{eq:3.11}, we get 
\begin{equation*}
E_\psi[W](s) \leq L_3\left(1+ \left(E[W](0)\right)^{1/2}\right)+L_4
\end{equation*}
for $s\geq1$. For $s\leq1$ we have 
\begin{equation*}
E_\psi[W](s) \leq L_3\left(1+ \left(E[W](0)\right)^{1/2}\right)+ E_\psi[W](0)
\end{equation*}
in the same way. 
Therefore inequality \eqref{eq:3.1} follows,  
and the proof is complete. 
$\Box$

\noindent{\bf Proof of Proposition \ref{Proposition:3.2}.}
By \eqref{eq:3.10}, \eqref{eq:2.2} and Jensen's inequality we have
\begin{equation*}
\begin{split}
\frac{1}{2}&\frac{d}{ds}\int_{{\bf R}^N}  \psi^2 |W|^2 \rho \,dy  
 \geq -2E_\psi [W](s) + c_G (p-1)\int_{{\bf R}^N} \psi^2 |W|^{p+1}\rho \,dy -L_5 \\
& \geq  -2E_\psi[W](s) + c_2 \left( \int _{{\bf R}^N} \psi^2 |W|^2 \rho \,dy\right)^{\frac{p+1}{2}} -L_5, 
\end{split}
\end{equation*}
where $c_2=(4\pi)^{-\frac{N(p-1)}{4}}(p-1)c_G \| \psi\|_\infty ^{-(p-1)}$.
Set 
\begin{equation*}
T_1:=\int_0^\infty \frac{d\tau}{1+2c_2\tau^{\frac{p+1}{2}}},
\end{equation*}
and we show $E_\psi[W](s)\geq -L_2$ for $L_2= L_3(T_1+E[W](0)^{1/2})+L_5/2 +1/4$. If not, then there exists some point $s_2\in (0,\infty)$ such that $E_\psi[W](s_2)< -L_2$.

For $s\in[0,T_1]$, as in the proof of Proposition~\ref{Proposition:3.1}, by \eqref{eq:3.8} we have 
\begin{equation*}
\begin{split}
E_\psi[W](s_2+s)
& <-L_2 +L_3 \int_{s_2}^{s_2+s} (1+\|W_s(\tau);L_\rho^2({\bf R}^N)\|) \,d\tau\\
& \leq -L_2 + L_3T_1+ L_3 E[W](0)^{1/2}\\
& = -\frac{L_5}{2}- \frac{1}{4}.
\end{split}
\end{equation*}
Therefore, for $s\in[s_2,s_2+T_1]$, we get
\begin{equation*}
\frac{d}{ds}\int_{{\bf R}^N}  \psi^2 |W|^2 \rho \,dy  
 \geq  1 + 2c_2 \left( \int _{{\bf R}^N} |W|^2 \rho \,dy\right)^{\frac{p+1}{2}}.
\end{equation*}
Since the solution to $f' =1+2c_2f^{(p+1)/2}$ and $f(0)=0$ blows up at $T_1$, this is a contradiction. Therefore we obtain \eqref{eq:3.2}, and the proof is complete.
$\Box$ 

\section{Proof of Theorem \ref{Theorem:1.1}.}
In this section we obtain \eqref{eq:1.11} by a bootstrap argument in Proposition~\ref{Proposition:4.1} and prove Theorem~\ref{Theorem:1.1}. Fix $\phi \in C^\infty ([0,\infty))$ such that $\phi (t)=1$ for $t\leq 1$ and $\phi(t)=0$ for  $t\geq 2$. For $R>0$, set $\varphi$  
\begin{equation*}
\varphi(x) = \phi\left(\frac{|x|}{R}\right), \quad x\in {\bf R}^N.
\end{equation*}
Let $B_R$ be an open ball of radius $R$ centered at the origin of ${\bf R}^N$.
We give the key estimate for the proof of Theorem~\ref{Theorem:1.1}. 

\begin{proposition}
\label{Proposition:4.1} 
Assume that $1<p<p_S$. Let $W$ satisfy \eqref{eq:1.6} with \eqref{eq:1.8}.  For $q\geq 2$ and $R>0$, there exists  $C_{q,R}>0$ such that 
\begin{equation}
\label{eq:4.1}
\sup _{s\geq 0} \int _s^{s+1} \|W; W_\rho ^{1,2} (B_R)\| ^{2q} \,d\tau \leq C_{q,R}.
\end{equation}
\end{proposition}
First we prove Proposition~\ref{Proposition:4.1} for the case $q=2$. 

\noindent{\bf Proof of Proposition~\ref{Proposition:4.1} for q=2.}
For any $R>0$ and $s\geq0$, we see by \eqref{eq:2.8} and \eqref{eq:2.11} that 
\begin{equation*}
\begin{split}
\int _s^{s+1} \|W; W_\rho ^{1,2} (B_R)\| ^{4} \,d\tau 
& \leq 2K_3 ^2  \left (1+ \int _s^{s+1} \|W_s; L_\rho ^{2} ({\bf R}^N)\| ^{2} \,d\tau \right)\\
& \leq 2K_3 ^2 (1+E[W](0)).
\end{split}
\end{equation*}
Therefore Proposition~\ref{Proposition:4.1} follows for $q=2$.  
$\Box$

Hereafter we assume $q\geq2$. 
As in \cite{GMS2}, we introduce several constants $p_1$, $\overline{q}$, $\lambda_q$, $\lambda$, $\theta$ and $\alpha$, satisfying
\begin{equation}
\label{eq:4.2}
\left\{
\begin{array}{ll}
& \displaystyle p_1 =1+ \frac{1}{p}, \quad q<\overline q<q+\frac{1}{p+1}, \quad  \lambda _q = p+1 - \frac{p-1}{q+1},\vspace{3pt}\\
& \displaystyle 2<\lambda<\lambda _q, \quad \theta = \frac{(p+1)(\lambda -2)}{(p-1)\lambda}, \quad  \alpha = \frac{2}{(1-\theta)\overline q}. 
\end{array}
\right.
\end{equation}
 
\begin{proposition}
\label{Proposition:4.2}
Assume that $q\geq2$. Let $p_1$, $\overline{q}$, $\lambda_q$, $\lambda$, $\theta$ and $\alpha$ be as in \eqref{eq:4.2}. Then there exists $\lambda_0 \in (2, \lambda_q)$ depending only on $p$ and $q$ such that 
\begin{equation}
\label{eq:4.3}
\alpha > 1, 
\end{equation}
\begin{equation}
\label{eq:4.4}
1<\frac{\theta \overline q \alpha'}{p_1}<q, 
\end{equation}
where $\alpha'$ is the  H\"{o}lder conjugate   of $\alpha$, i.e.,  $1/\alpha + 1/\alpha' = 1$. 
\end{proposition}
{\bf Proof.}
The condition $\alpha >1$ is equivalent to 
\begin{equation*}
\overline q < \frac{2}{1-\theta}=\frac{(p-1)\lambda}{(p+1)-\lambda}.
\end{equation*}
Since the right hand side is monotone increasing for $\lambda \in (2, \lambda_q)$ and  
\begin{equation*}
\frac{(p-1)\lambda_q}{(p+1)-\lambda_q}=(p+1)\left (q+\frac{2}{p+1}\right)>q+\frac{1}{p+1}, 
\end{equation*}
we can choose $\lambda'=\lambda'(p,q)$ such that \eqref{eq:4.3} holds true for $\lambda>\lambda '$. 

We prove \eqref{eq:4.4}. Since
\begin{equation*}
\theta\overline q \alpha' = \frac{2\theta \overline q}{2-(1-\theta)\overline q}, 
\end{equation*}
the condition $p_1<\theta \overline{q}\alpha'$ is equivalent to
\begin{equation*}
\overline q > \frac{2p_1}{(2-p_1)\theta + p_1}. 
\end{equation*}
This together with 
\begin{equation*}
\frac{2p_1}{(2-p_1)\theta + p_1}<2<\overline q 
\end{equation*}
implies that $p_1<\theta \overline{q}\alpha'$ for any $\lambda>\lambda '$. 
 On the other hand, the condition $\theta \overline{q}\alpha'<p_1 q$ is equivalent to  
\begin{equation*}
\overline q <\frac{2p_1 q}{(2-p_1 q)\theta + p_1 q}. 
\end{equation*}
Since the right hand side is monotone increasing for $\lambda \in (2, \lambda_q)$ and
\begin{equation*}
\frac{2p_1 q}{(2-p_1 q)\theta_{\lambda_q} + p_1 q} = q + \frac{2}{p+1} > q+\frac{1}{p+1}, 
\end{equation*}
we can choose $\lambda''=\lambda''(p,q)$ so that \eqref{eq:4.4} holds for $\lambda>\lambda ''$. 
\eqref{eq:4.3} and \eqref{eq:4.4} holds for  $\lambda_0 > \max\{\lambda', \lambda''\}$. Therefore the proof is complete.
$\Box$

Hereafter we fix $\lambda =\lambda_0$ so that relations \eqref{eq:4.2},  \eqref{eq:4.3} and \eqref{eq:4.4} hold.  To prove Proposition~\ref{Proposition:4.1}, we prepare the following four lemmas.

\begin{lemma}
\label{Lemma:4.1}
Assume that $1<p<p_S$. Let $W$ satisfy \eqref{eq:1.6} with \eqref{eq:1.8}.  Assume that inequality \eqref{eq:4.1} holds with $B_R$ replaced by $B_{2R}$ for some $q\geq 2$. Let $p_1$, $\overline{q}$, $\lambda_q$, $\lambda$, $\theta$ and $\alpha$ be as in \eqref{eq:4.2}. Then there exists $J_1>0$ such that 
\begin{equation}
\label{eq:4.5}
\begin{split}
\int _\sigma ^{s+1} \|\varphi W_s (\tau);& L_\rho ^{p_1} (B_{2R})\| ^{\theta \overline q \alpha' } \,d\tau
 \leq J_1 \left(1+ \int _\sigma ^{s+1} \| |W|^{p}; L_\rho ^{p_1} (B_{2R})\| ^{\theta \overline q \alpha' } \,d\tau \right)
\end{split}
\end{equation}
for all $s\ge 1$ with some $\sigma \in [s-1/4, s]$. 
\end{lemma} 
{\bf Proof.}
Let $\sigma \in [s-1/4, s]$. 
We apply the maximal regularity theorem for the heat equation \cite[Lemma~6.5]{GMS2} to 
\begin{equation*}
\label{eq:4}
\begin{split}
(\varphi w_{is}) & - \Delta (\varphi w_i) 
 = -2 \nabla \varphi \cdot \nabla w_i -\Delta \varphi w_i -\frac{\varphi}{2}y\cdot \nabla w_i -\beta \varphi w_i + \varphi F_i(W)=:f_i
\end{split}
\end{equation*}
in $B_{2R}\times (\sigma, s+1)$ to get  
\begin{equation}
\label{eq:4.6}
\begin{split}
& \int _\sigma ^{s+1} 
 \|\varphi W_s; L_\rho ^{p_1} (B_{2R})\| ^{\theta \overline q \alpha' } \,d\tau \\
& \leq C_* \left( \int_\sigma ^{s+1} \|f;L_\rho ^{p_1} (B_{2R})\|^{\theta \overline q\alpha'}\,d\tau + \|W(\sigma); C^2(\overline {B_{2R}})\|^{\theta \overline q\alpha'}\right),  
\end{split}
\end{equation}
where $f=(f_1,\dots,f_M)$ and  $C_*$ is a positive constant. Let $c_3=(4\pi)^{\frac{N(p-1)}{4(p+1)}}$. We see by H\"older's inequality and \eqref{eq:2.9} that 
\begin{equation}
\label{eq:4.7}
\begin{split}
\|\Delta \varphi W ; L_\rho ^{p_1}(B_{2R})\|
& \leq \left(\int_{B_{2R}}|W|^2 \rho \,dy\right)^{\frac{1}{2}}\left(\int_{B_{2R}}|\Delta \varphi|^{2(p+1)/(p-1)}\rho \,dy\right)^{\frac{p-1}{2(p+1)}} \\
& \leq c_3 K_1\|\Delta \varphi \|_{\infty}.
\end{split}
\end{equation}
In the same way, we obtain 
\begin{align}
\label{eq:4.8}
& \|\beta \varphi W; L_\rho ^{p_1}(B_{2R})\|\leq \beta c_3 K_1, \\
\label{eq:4.9}
& \||\nabla \varphi||\nabla W| ; L_\rho ^{p_1}(B_{2R})\|\leq c_3 \|\nabla \varphi \|_\infty \| \nabla W ; L_\rho ^2 (B_{2R})\|, \\
\label{eq:4.10}
& \left\|\frac{\varphi}{2} |y| |\nabla W |; L_\rho ^{p_1}(B_{2R})\right\|\leq c_3 R \| \nabla W ; L_\rho ^2 (B_{2R})\|.
\end{align}
Since $\theta \overline q \alpha'<2q$, by \eqref{eq:4.9} we see that 
\begin{equation}
\label{eq:4.11}
\begin{split}
& \left (\int _\sigma ^{s+1} 
 \|2|\nabla \varphi|| \nabla W|; L_\rho ^{p_1} (B_{2R})\| ^{\theta \overline q \alpha' } \,d\tau \right)^{1/\theta \overline q \alpha' }\\
& \leq 4c_3 \|\nabla \varphi \|_\infty \left (\int _\sigma ^{s+1} \| \nabla W; L_\rho ^{2} (B_{2R})\| ^{2q} \,d\tau \right)^{1/2q}\\  
& \leq 8c_3 \|\nabla \varphi \|_\infty C_{q,2R} ^{1/2q}.
\end{split}
\end{equation}
In the same way, by \eqref{eq:4.10} we obtain 
\begin{equation}
\label{eq:4.12}
\left (\int _\sigma ^{s+1} \left\|\frac{\varphi}{2}|y| |\nabla W|; L_\rho ^{p_1} (B_{2R})\right\| ^{\theta \overline q \alpha' } \,d\tau \right)^{1/\theta \overline q \alpha' }\leq  4c_3  R C_{q,2R} ^{1/2q}. 
\end{equation}
By \eqref{eq:2.2} we  see that
\begin{equation}
\label{eq:4.13}
\begin{split}
& \left (\int_\sigma ^{s+1} \| \varphi F(W);  L_\rho ^{p_1}(B_{2R})\| ^{\theta \overline q \alpha'}\,d\tau \right)^{1/\theta \overline q \alpha'}\\
& \leq C_F \left(\int_\sigma ^{s+1} \| |W|^p; L_\rho ^{p_1}(B_{2R})\|^{\theta \overline q \alpha'} \,d\tau \right)^{1/\theta \overline q \alpha'}.
\end{split}
\end{equation}
Thus by \eqref{eq:4.7}, \eqref{eq:4.8}, \eqref{eq:4.11}, \eqref{eq:4.12} and \eqref{eq:4.13}  we obtain
\begin{equation}
\label{eq:4.14}
\begin{split}
& \left (\int _\sigma ^{s+1} 
 \|f; L_\rho ^{p_1} (B_{2R})\|^{\theta \overline q \alpha'}\,d\tau \right)^{1/\theta \overline q \alpha'}\\
& \leq 8c_3  \|\nabla \varphi \|_\infty C_{q,2R} ^{1/2q}+  4c_3  R C_{q,2R}^{1/2q}+ 2c_3 K_1\|\Delta \varphi \|_\infty \\
& + 2\beta c_3 K_1 +   C_F \left (\int_\sigma ^{s+1} \| |W|^p; L_\rho ^{p_1}(B_{2R})\|^{\theta \overline q \alpha'} \,d\tau \right)^{1/\theta \overline q \alpha'}. 
\end{split}
\end{equation}
As in \cite[Lemma~6.9 and Remark~6.10]{GMS2}, by \eqref{eq:2.8} and \eqref{eq:2.11} we can choose $\sigma$ such that $\|W(\sigma); C^2(\overline {B_{2R}})\|\leq Z$ with a constant $Z$ independent of $s$ and $a$. 
This together with \eqref{eq:4.6} and \eqref{eq:4.14} implies \eqref{eq:4.5},  and the proof is complete.
$\Box$

\begin{lemma}
\label{Lemma:4.2}
Let $W$ satisfy \eqref{eq:1.6} with \eqref{eq:1.8}. For any $R>0$, then there exists  $J_2>0$ such that
\begin{equation}
\label{eq:4.15}
\|W(s); W_\rho ^{1,2} (B_R)\| ^{2}\leq J_2(1+\| \varphi^2|W(s)||W_s(s)|;L^1_\rho (B_{2R})\|) 
\end{equation}
for $s\ge0$. 
\end{lemma}
{\bf Proof.}
It follows from \eqref{eq:1.10} that 
\begin{equation*}
\begin{split}
-2E_\varphi[W](s) +(p-1)& \int_{{\bf R}^N} \varphi^2G(W)\rho \,dy \\
& = \frac{p-1}{2} \int _{{\bf R}^N} \varphi^2 (|\nabla W |^2+\beta |W|^2)\rho \,dy-(p+1)E_\varphi[W](s).
\end{split}
\end{equation*}
This together with \eqref{eq:3.10} implies that 
\begin{equation*}
\begin{split}
&\|W; W_\rho ^{1,2}  (B_R)\| ^{2} 
 \leq \int _{{\bf R}^N} \varphi^2 (|\nabla W |^2+\beta |W|^2)\rho \,dy\\
& \leq \frac{2}{p-1}((p+1)L_1 + L_5 +\| \varphi^2|W||W_s|;L^1_\rho (B_{2R})\|)\\
& \leq J_3(1+\| \varphi^2|W||W_s|;L^1_\rho (B_{2R})\|).
\end{split}
\end{equation*}
Thus inequality \eqref{eq:4.15} follows. 
$\Box$

\begin{lemma}
\label{Lemma:4.3}
Let $W$ satisfy \eqref{eq:1.6} with \eqref{eq:1.8}. Assume that  for some $q\geq 2$, \eqref{eq:4.1} holds with $B_R$ replaced by $B_{2R}$. Let $\lambda_q$ be as in \eqref{eq:4.2}. Then there exist  $C_{q,R}'>0$ and  $C_{q,R}''>0$ such that 
\begin{align}
\label{eq:4.16}
& \sup _{s\geq0} \int _s^{s+1} \|W; L_\rho ^{p+1} (B_{R} )\| ^{(p+1)q} \,d\tau \leq C_{q,R}', \\
\label{eq:4.17}
& \sup _{s\geq 0} \|W; L_\rho ^{\lambda } (B_{R})\| \leq C_{q,R}'' \quad {\rm for}\,\,{\rm all}\,\, \lambda < \lambda _q.  
\end{align}
\end{lemma}
{\bf Proof.}
By \eqref{eq:2.2} and \eqref{eq:3.2} we see that for any $R>0$, 
\begin{equation}
\label{eq:4.18}
\begin{split}
\|W; L_\rho ^{p+1} (B_R)\| ^{(p+1)}
& \leq \frac{1}{c_G} \int _{{\bf R}^N} \varphi^2 G(W)\rho \,dy \\
& \leq \frac{1}{c_G} \left(L_2 + \|W; W_\rho ^{1,2} (B_{2R})\| ^{2}\right). 
\end{split}
\end{equation}
Thus we get 
\begin{equation*}
\begin{split}
\int_s ^{s+1} \|W; L_\rho ^{p+1} (B_R)\| ^{(p+1)q}\,d\tau
& \leq \frac{2^{q-1}}{c_G^q}\int_s ^{s+1}  (L_2^q + \|W; W_\rho ^{1,2}(B_{2R})\|^{2q})\,d\tau\\
& \leq \frac{2^{q-1}}{c_G^q} (L_2^q + C_{q, 2R}).
\end{split}
\end{equation*}
Therefore  inequality \eqref{eq:4.16} follows. 
Furthermore, by \eqref{eq:4.16} and \eqref{eq:2.8} we obtain 
\begin{equation*}
\begin{split}
\sup_{s\geq0} \int_s^{s+1}( \|W;L ^{p+1} (B_{R})\|^{(p+1)q}& + \|W_s;L ^2(B_{R}) \|^2)\,d\tau \\
& \leq \exp\left(\frac{qR^2}{4}\right)\left(C_{q,R}' + E[W](0)\right).
\end{split}
\end{equation*}

 This together with the interpolation theorem \cite[Lemma~A.1]{GMS2} (See also \cite{CL}) implies \eqref{eq:4.17}. The proof is complete. 
$\Box$
\begin{lemma}
\label{Lemma:4.4}
Assume that $1<p<p_S$. Let $W$ satisfy \eqref{eq:1.6} with \eqref{eq:1.8}.  Assume that  for some $q\geq 2$, \eqref{eq:4.1} holds with $B_R$ replaced by $B_{4R}$. Let $p_1$, $\overline{q}$, $\lambda_q$, $\lambda$, $\theta$ and $\alpha$ as in \eqref{eq:4.2}. Then there exists $J_3>0$ such that 
\begin{equation}
\label{eq:4.19}
\begin{split}
\int _\sigma ^{s+1} \|W;W_\rho ^{1,2} (B_R)\| ^{2\overline q} \,d\tau 
\leq J_3 \left (1+ \int _\sigma ^{s+1} \|W; W_\rho ^{1,2} (B_{4R})\| ^{2\theta \overline q \alpha' /p_1} \,d\tau \right)
\end{split}
\end{equation}
for all $s\geq 1$ with some $\sigma \in [s-1/4, s]$. 
\end{lemma}
{\bf Proof.}
By \eqref{eq:4.15}, \eqref{eq:4.17} and H\"older's inequality we have  
\begin{equation*}
\begin{split}
\|W; W_\rho ^{1,2} (B_{R})\| ^{2}
& \leq J_2 (1+\| \varphi W;L_\rho ^{\lambda} (B_{2R})\| \| \varphi W_s;
L_\rho ^{\lambda'} (B_{2R})\|)\\
& \leq J_2 (1+C_{q,2R} ''\| \varphi W_s;L_\rho ^{p_1} (B_{2R})\|^\theta  \| \varphi W_s;
L_\rho ^{2} (B_{2R})\|^{1-\theta}).
\end{split}
\end{equation*}
Then we obtain 
\begin{equation}
\label{eq:4.20}
\begin{split}
& \int _\sigma  ^{s+1}
  \|W; W_\rho ^{1,2} (B_{R})\| ^{2\overline q} \,d\tau \\
& \leq  2^{\overline q-1}J_2^{\overline q} \left (2+ C_{q,2R}''^{\overline q}\int _\sigma ^{s+1} \| \varphi W_s;L_\rho ^{p_1} (B_{2R})\|^{\theta \overline q}  \| \varphi W_s;L_\rho ^{2} (B_{2R})\|^{(1-\theta)\overline q} \,d\tau \right).
\end{split}
\end{equation}
By \eqref{eq:2.8}, \eqref{eq:4.2} and H\"older's inequality  we see that 
\begin{equation}
\label{eq:4.21}
\begin{split}
& \int _\sigma ^{s+1} 
 \| \varphi W_s;L_\rho ^{p_1} (B_{2R})\|^{\theta \overline q}  \| \varphi W_s;L_\rho ^{2} (B_{2R})\|^{(1-\theta)\overline q} \,d\tau \\
& \leq (E[W](0))^{1/\alpha}\left (\int _\sigma ^{s+1} \| \varphi W_s;L_\rho ^{p_1} (B_{2R})\|^{\theta \overline q\alpha'} \,d\tau \right)^{1/\alpha'}.
\end{split}
\end{equation} 
By \eqref{eq:4.18} we have 
\begin{equation}
\label{eq:4.22}
\|W;L_\rho ^{p+1}(B_R)\|^{p+1}
\leq J_4(1+\|W;W_\rho ^{1, 2}(B_{2R})\|^2 )
\end{equation}
for $J_4:= c_G^{-1} \max\{L_2,1\}$.
Furthermore, by \eqref{eq:4.5} and \eqref{eq:4.22} we have
\begin{equation}
\label{eq:4.23}
\begin{split}
& \left (\int _\sigma ^{s+1} 
 \| \varphi W_s;L_\rho ^{p_1} (B_{2R})\|^{\theta \overline q\alpha'} \,d\tau \right)^{1/\alpha'}\\
& \leq J_1 ^{1/\alpha'}\left (1+ \int _\sigma ^{s+1} \| |W|^{p}; L_\rho ^{p_1} (B_{2R})\| ^{\theta \overline q \alpha' } \,d\tau \right)\\
& = J_1 ^{1/\alpha'} \left (1+ \int _\sigma ^{s+1} \| W; L_\rho ^{p+1} (B_{2R})\| ^{p \theta \overline q \alpha' } \,d\tau \right)\\
& \leq J_1 ^{1/\alpha'}  \left(1+ 2^{\theta \overline q \alpha'/p_1 -1}J_4^{\theta \overline q \alpha'/p_1}\left(2+ \int _\sigma^{s+1} \| W; W_\rho ^{1,2} (B_{4R})\| ^{2\theta \overline q \alpha' /p_1} \,d\tau \right)\right). 
\end{split}
\end{equation}
Therefore inequality \eqref{eq:4.19} follows from \eqref{eq:4.20}, \eqref{eq:4.21} and \eqref{eq:4.23}. 
$\Box$

\noindent{\bf Proof of Proposition~\ref{Proposition:4.1}.}
We have already proved Proposition~\ref{Proposition:4.1} for $q=2$. We prove Proposition~\ref{Proposition:4.1} for $\overline{q}>2$ by a bootstrap argument starting with $q=2$. If \eqref{eq:4.1} holds for some $q\geq 2$ with some $4R>0$, then by \eqref{eq:4.4} and \eqref{eq:4.19} we have 
\begin{equation*}
\begin{split}
\int _s ^{s+1} \|W; W_\rho ^{1,2} (B_R)\| ^{2\overline q} \,d\tau 
& \leq J_3 \left (1+ 2\left (\int _\sigma ^{s+1} \|W; W_\rho ^{1,2} (B_{4R})\| ^{2q} \,d\tau \right)^{\theta \overline q \alpha' /qp_1}\right)\\
& \leq J_3 \left(1+ 2 (2C_{q,4R})^{\theta \overline q \alpha' /qp_1} \right)
\end{split}
\end{equation*}
for $\overline q < q+1/(p+1)$ and $s\geq 1$. On the other hand, we have 
\begin{equation*}
\sup _{0\leq s\leq 1} \int _s^{s+1} \|W; W_\rho ^{1,2} (B_{R})\| ^{2\overline{q}} \,d\tau \leq C\sup_{0\leq t\leq 1-e^{-2}}\|U(t);BC^1({\bf R}^N)\|^{2\overline{q}}.
\end{equation*}
Therefore we see that \eqref{eq:4.1} holds for $\overline q < q+1/(p+1)$ with $R$.  
Let $q_1>2$,  $R_1>0$ and $m=[(p+1)({q_1}-2)]+1$. We repeat the estimates $m$ times starting with $q=2$ and $R=4^mR_1$ to obtain 
\begin{equation*}
\sup _{s\geq 0} \int _s^{s+1} \|W; W_\rho ^{1,2} (B_{R_1})\| ^{2q_1} \,d\tau \leq C_{q_1,R_1}.
\end{equation*}
Therefore the proof of \ref{Proposition:4.1} is complete. 
$\Box$

\noindent{\bf Proof of Theorem~\ref{Theorem:1.1}.}
We use the interior regularity theorem for linear parabolic equation. 
Let $w_i^{\pm}:=\max\{\pm w_i,0\}$.
Since $W=(w_1,\dots,w_M)$ satisfies
\begin{equation*}
\partial _s w_i = \Delta w_i- \frac{1}{2}y\cdot \nabla w_i -\beta w_i + \sum_{j=1} ^M \beta_{ij} |w_i|^{r-1} |w_j|^{r+1}w_i\quad \mbox{in} \quad {\bf R}^N\times (0,\infty)
\end{equation*}
for $i=1,\dots, M$, $\{w_i^{\pm}\}_{i=1} ^M$ satisfies
\begin{equation*}
\partial _s w_i^{\pm} \leq \Delta w_i^{\pm}- \frac{1}{2}y\cdot \nabla w_i^{\pm} -\beta w_i^{\pm} + \sum_{j=1} ^M \beta_{ij} |w_i|^r|w_i|^{r+1}\quad \mbox{in} \quad {\bf R}^N\times (0,\infty)
\end{equation*}
in a weak sense. See e.g.,  \cite[Chapter 1]{D}. Therefore 
\begin{equation*}
w:= \sum_{i=1}^M |w_i|
\end{equation*}
satisfies
\begin{equation}
\label{eq:4.24}
\partial _s w \leq \Delta w- \frac{1}{2}y\cdot \nabla w -\beta w + \sum_{i,j=1} ^M \beta_{ij} w^{p}
\end{equation}
in ${\bf R}^N\times (0,\infty)$ in a weak sense.
We choose $\alpha$,  $\beta$ and $q$ so that $1/\beta + N/2\alpha <1$, $\alpha \geq 1$ and $\alpha (p-1)<\lambda_1 (q)$. This is possible for sufficiently large $q$ provided that $1<p<p_S$.
By \eqref{eq:2.9} and \eqref{eq:4.17} we can apply  an interior regularity theorem for a linear parabolic equation \cite[Lemma~A.2]{GMS2} (see also \cite{LSU}) for \eqref{eq:4.24} in $B_R\times (s,s+1)$ for $s\geq0$ to find  $C>0$ such that 
\begin{equation*}
\|w;L^\infty (B_{R/2})\|\leq C \quad {\rm for}\,\,s> \frac{1}{2}.
\end{equation*}
This together with \eqref{eq:1.8} implies 
\begin{equation*}
\|W;L^\infty (B_{R/2})\|\leq C \quad {\rm for}\,\,s\geq 0.
\end{equation*}
Since the constant $C$ is independent of $a$, we have
\begin{equation*}
|U(x,t)|\leq C(1-t)^{-\beta} \quad {\rm for}\,\,0<t<1,
\end{equation*}
and the proof of Theorem~\ref{Theorem:1.1} is complete.   
$\Box$
\vspace{7pt}

\noindent
{\bf Acknowledgment.} 
The author was supported in part by JSPS KAKENHI Grant Number 20J11261.


\bibliographystyle{amsplain}

\bigskip
\noindent Addresses:

\smallskip
\noindent E. Z.: Graduate School of Mathematical Sciences, The University of Tokyo 3-8-1 Komaba, Meguro-ku, Tokyo 153-8914, Japan\\
\noindent 
E-mail: {\tt erbol@ms.u-tokyo.ac.jp}\\
\end{document}